\documentclass[twocolumn]{autart}    

\makeatletter
\renewcommand\@biblabel[1]{}
\renewenvironment{thebibliography}[1]
     {\section*{\refname}%
      \@mkboth{\MakeUppercase\refname}{\MakeUppercase\refname}%
      \list{\@biblabel{\@arabic\c@enumiv}}%
           {\settowidth\labelwidth{\@biblabel{#1}}%
            \leftmargin\labelwidth
            \advance\leftmargin\labelsep
             \advance\leftmargin by 0.5cm%
            \itemindent -0.6cm
            \usecounter{enumiv}%
            \let\p@enumiv\@empty
            \renewcommand\theenumiv{\@arabic\c@enumiv}}%
      \sloppy
      \clubpenalty4000
      \@clubpenalty \clubpenalty
      \widowpenalty4000%
      \sfcode`\.\@m}
     {\def\@noitemerr
     {\@latex@warning{Empty `thebibliography' environment}}%
     \endlist}
\makeatother

\usepackage{graphicx}          
\usepackage{mathrsfs}
\usepackage{amsmath,amssymb,color}
\newenvironment{shrinkeq}[1]
{ \bgroup
  \addtolength\abovedisplayshortskip{#1}
  \addtolength\abovedisplayskip{#1}
  \addtolength\belowdisplayshortskip{#1}
  \addtolength\belowdisplayskip{#1}}
{\egroup\ignorespacesafterend}

\begin{document}
\begin{frontmatter}

\title{Decomposition with respect to outputs \\ for Boolean control
networks.\thanksref{footnoteinfo}}

\thanks[footnoteinfo]{The research is supported in part by National Natural Science
Foundation (NNSF) of China under Grants 11271194, 61273115,
61074115 and 10701042.}

\author{Yunlei Zou}\ead{zouyl0903@163.com},               
\author{Jiandong Zhu}\ead{jiandongzhu@njnu.edu.cn}

\address{School of Mathematical Sciences, Nanjing Normal
University, Nanjing, 210023, PR~China}

\begin{abstract}                
This paper investigates the problem of decomposition with respect
to outputs for Boolean control networks (BCNs). First, with the
linear expression of BCNs and the matrix semi-tensor product,
some algebraic equivalent conditions for the decomposition are
obtained. Second, a necessary and sufficient graphical condition
for the decomposition with respect to outputs is given. Third,
an effective method is proposed to reduce the computational burden in the
realization of the decomposition. Finally, some examples are
addressed to validate the effectiveness of the proposed method.
\end{abstract}


\end{frontmatter}

\section{Introduction}

Boolean networks (BNs) are a kind of discrete dynamical systems
described by logical variables and logical functions. They are
first proposed by Kauffman (1969) to describe and analyze cell
regulation (Albert \& Othmer, 2003; Faure, Naldi, Chaouiya, \&
Thieffry, 2006). BNs with additional inputs and outputs are
usually called Boolean control networks (BCNs), and they have been
paid great attention by biologists and control theory scientists
(Akutsu, Hayashida, Ching, \& Ng, 2007; Cheng, 2009;). In recent
years, a semi-tensor product method of BNs and BCNs was developed
by Cheng and his collaborators (Cheng, 2009; Cheng \& Qi, 2009;
Cheng \& Qi, 2010a). For a complete control theory framework of
the BCNs, we refer to the book written by Cheng, Qi, and Li
(2011). Based on the linear algebraic expression of BCNs, many
classical problems in control theory can be generalized to BCNs
such as controllability, observability, stabilization, disturbance
decoupling and optimal control (Cheng \& Qi, 2009;
Laschov \& Margaliot, 2012; Laschov, Margaliot, \& Even, 2013; Fornasini \& Valcher, 2013a;
Fornasini \& Valcher, 2013b; Cheng, 2011; Zhao, Li \& Cheng, 2011). Moreover, some
of these results have be further extended to different kinds of
BCNs (Li, Yang, \& Chu, 2014; Li \& Sun, 2012; Feng, Yao, \& Cui,
2013; Li \& Wang, 2012).
\par
System decomposition is an important issue in the traditional
linear control system theory, which makes clear the observable and
controllable components of a system. An interesting topic is
ascertaining whether there exists a similar system decomposition for
BCNs. As a matter of fact, this topic has been discussed by Cheng, Li
and Qi (2010), where the controllable normal form, the observable
normal form and the Kalman decomposition form of BCNs are
presented with the state-space analysis method (Cheng \& Qi,
2010b). The state space of an $n$-dimensional BCN is
defined as the set of all logical functions with respect to the
$n$ logical variables. Since the cardinal of the state space of an
$n$-dimensional BCN is $2^{2^n}$, generally speaking, the
computational burden of state space method is very heavy for large
$n$ although it displays great advantage in theoretical analysis.
Furthermore, some additional regularity assumptions are imposed on
the system model for achieving those decompositions. \par In our
recent paper Zou \& Zhu (2014), we have proposed the {\it
decomposition with respect to inputs} for BCNs, which is just the
controllable normal form proposed by Cheng, Li and Qi (2010) under
the regularity assumption on the largest uncontrollable subspace.
But we did not use any concepts of state spaces. A main contribution of Zou
\& Zhu (2014) lies in that the system decomposition can be realized
without the regularity condition on the largest uncontrollable
subspace. Moreover, a graphical condition and a constructive
algorithm for the decomposition with respect to inputs are
proposed. Then, a natural idea is reconsidering the observable
normal form and the Kalman decomposition form for BCNs.
\par
In this paper, we consider the decomposition with respect to
outputs in the framework of the linear algebraic representation of
BCNs. First, some necessary and sufficient algebraic conditions
for the decomposition with respect to outputs are proposed. Next,
an equivalent graphical condition for the decomposition with
respect to outputs is derived. As the graphical
condition is satisfied, a constructive coordinate transformation
to realize the decomposition is obtained. Then, a method is
proposed to reduce the computational burden in realizing the
decomposition with respect to outputs. Finally, some
examples are analyzed with our method. \par We organize this paper
as follows. In section 2, some preliminaries and problem statement
are provided. In section 3, some algebraic conditions for the
decomposition with respect to outputs are presented. In section 4,
a graphical condition is obtained. In Section 5, a computation method is given. Finally, Section 6 gives a brief conclusion.
\vspace{-3mm}

\section{Preliminaries and Problem Statement}

Let $\mathcal{D}=\{\mbox{True}=1,
\mbox{False}=0\}$.
Consider a BCN described by the logical equations
\vspace{-3mm}
\begin{eqnarray}
  x_1(t + 1) &=& f_1(x_1(t),\cdots, x_n(t), u_1(t),\cdots, u_m(t)), \nonumber\\
    &\vdots&\nonumber\\
    x_n(t + 1) &=& f_n(x_1(t),\cdots, x_n(t), u_1(t),\cdots, u_m(t)), \nonumber\\
    \label{eq1}y_1(t) &=& h_1(x_1(t),\cdots, x_n(t)), \\
    &\vdots& \nonumber\\
    y_p(t) &=& h_p(x_1(t),\cdots, x_n(t)),\nonumber
\end{eqnarray}
where the state variables $x_i$, the output variables $y_i$ and the controls $u_i$ take values in $\mathcal{D}$,
$f_i: \mathcal{D}^{n+m}\rightarrow \mathcal{D}$ and $h_i: \mathcal{D}^{n}\rightarrow \mathcal{D}$ are logical functions.\par
Consider the logical mapping $G:\mathcal{D}^n\rightarrow\mathcal{D}^n$ defined by
\begin{eqnarray}
\label{eq2}
\hspace{2cm}z_1&=&g_1(x_1,x_2,\cdots, x_n),\nonumber \\
\hspace{2cm}z_2&=&g_2(x_1,x_2,\cdots, x_n),\nonumber \\
\hspace{2cm}&\vdots &\\
\hspace{2cm}z_n&=&g_n(x_1,x_2,\cdots, x_n).\nonumber
\end{eqnarray}
If $G:\mathcal{D}^n\rightarrow\mathcal{D}^n$ is a bijection, it is called a logical coordinate transformation.
\par
We say that BCN (\ref{eq1}) is {\it decomposable with respect to outputs of order $n-s$}, if there exists a
logical coordinate transformation (\ref{eq2}), such that (\ref{eq1}) becomes
\vspace{-3mm}
\begin{eqnarray}
  z_1(t + 1) &=& \hat{f}_1(z_1(t),\cdots, z_s(t), u_1(t),\cdots, u_m(t)), \nonumber\\
    &\vdots&\nonumber\\
    z_s(t + 1) &=& \hat{f}_s(z_1(t),\cdots, z_s(t), u_1(t),\cdots, u_m(t)), \nonumber\\
    z_{s+1}(t + 1) &=& \hat{f}_{s+1}(z_1(t),\cdots, z_n(t), u_1(t),\cdots, u_m(t)), \nonumber\\
    \label{eq3}&\vdots&\\
    z_n(t + 1) &=& \hat{f}_n(z_1(t),\cdots, z_n(t), u_1(t),\cdots, u_m(t)), \nonumber\\
    y_1(t) &=& \hat{h}_1(z_1(t),\cdots, z_s(t)), \nonumber\\
    &\vdots& \nonumber\\
    y_p(t) &=& \hat{h}_p(z_1(t),\cdots, z_s(t)).\nonumber
\end{eqnarray}
BCN (\ref{eq3}) is called a {\it decomposition with respect to
outputs of order $n-s$}. A decomposition with respect to outputs
of the maximum order is called the {\it maximum decomposition
with respect to outputs}. We say that BCN (\ref{eq1}) is {\it
undecomposable with respect to outputs} if the order of the
maximum decomposition with respect to outputs is $0$.
\par
Denote the real number field by $\mathbf{R}$ and the set of all the $m\times n$ real matrices by $\mathbf{R}_{m\times n}$. Let Col$(A)$ be the set of all the columns of matrix $A$ and denote the $i$th column of $A$ by $\mbox{Col}_{i}(A)$. Set $\Delta_{k}=\{\delta^{i}_{k}|i=1,2,\cdots,k\}$, where
$\delta^{i}_{k}=\mbox{Col}_i(I_{k})$ with $I_{k}$ the $k\times k$ identity
matrix. For simplicity, we denote $\Delta:=\Delta_2=\{\delta_{2}^{1},\delta_{2}^{2}\}$. A matrix $L\in
\mathbf{R}_{m\times n}$ is called a logical matrix if
Col$(L)\subset \triangle_{m}$. Obviously, the logical matrix $L$ satisfies $\mathbf{1}_{m}^\mathrm{T}L=\mathbf{1}_{n}^\mathrm{T}$, where $\mathbf{1}_{n}$ denote the $n-$dimensional column vector whose entries are all equal to 1. The set of all the $m\times r$ logical
matrices is denoted by $\mathcal{L}_{m\times r}$. For simplicity,
we denote the logical matrix
$L=[\delta^{i_{1}}_{m},\delta^{i_{2}}_{m},\cdots,\delta^{i_{r}}_{m}]$
by $\delta_{m}[i_{1},i_{2},\cdots,i_{r}]$.
 Set
$A\in \mathbf{R}_{m\times n}$, $B\in \mathbf{R}_{p\times q}$, and
$\alpha=\mathrm{lcm}(n,p)$ be the least common multiple of $n$ and
$p$. The left semi-tensor product of $A$ and $B$ is defined as (Cheng, Qi, \& Li, 2011): $$
A\ltimes B=(A\otimes I_{\frac{\alpha}{n}})(B\otimes
I_{\frac{\alpha}{p}}),$$ where $\otimes$ is the Kronecker product. Obviously,
the left semi-tensor product is a generalization of the traditional matrix product.
So $A\ltimes B$ can be written as $AB$. To express the logical equations with
linear algebraic method, elements in $\mathcal{D}$ are identified with vectors
$\mbox{True}\sim\delta_{2}^{1}$ and $\mbox{False}\sim \delta_{2}^{2}$ (Cheng, Qi, \& Li, 2011).

\par

{\bf Proposition 1.} (Cheng \& Qi, 2010a) {\it Let $x_i$ and
$u_i$ take values in $\Delta$ and denote $x=\ltimes_{i=1}^nx_i$,
$y=\ltimes_{i=1}^py_i$, $u=\ltimes_{i=1}^{m}u_i$. Then the BCN
(\ref{eq1}) can be expressed in the algebraic form
\begin{shrinkeq}{-2ex}
\begin{equation}
\label{eq4}
  x(t+1)=Lu(t)x(t), \ \  y(t)=Hx(t),
\end{equation}
\end{shrinkeq}
where $L\in \mathcal{L}_{2^n\times 2^{n+m}}$ and $H\in \mathcal
{L}_{2^p\times 2^{n}}$.}

\par

Let $z=Tx$ be the algebraic form of logical coordinate
transformation (\ref{eq2}), where $T$ is a permutation matrix. Set
$z^{[1]}=\ltimes_{i=1}^sz_i$ and $z^{[2]}=\ltimes_{i=s+1}^nz_i$.
Then the decomposition form (\ref{eq3}) can be rewritten in the
algebraic form
\begin{equation}
\label{eq5}
\begin{array}{rcl}
z^{[1]}(t+1)& =& G_1u(t)z^{[1]}(t), \\
z^{[2]}(t+1)& =& G_2u(t)z(t),\\
y(t)& =& Mz^{[1]}(t),
\end{array}
\end{equation}
where $G_1\!\in\!\mathcal{L}_{2^s\times 2^{s+m}}$,
$G_2\!\in\! \mathcal{L}_{2^{n-s}\times2^{n+m}}$ and
$M\!\in\!\mathcal{L}_{2^{p}\times2^{s}}$.

\par

For the BCN (\ref{eq1}) with algebraic form (\ref{eq4}), the {\it
problem of decomposition with respect to outputs} is to find a
coordinate transformation matrix $T$ such that (\ref{eq4}) has the form
(\ref{eq5}). The {\it problem of maximum decomposition with
respect to outputs} is to find a coordinate transformation matrix $T$ such that
the decomposition with respect to outputs has the maximum order
$n-s$.
\par
In Cheng, Li and Qi (2010) , the observable normal form of a BCN
is proposed. Here, we rewrite the result in the algebraic form as follows.
\par
{\bf Proposition 2.} (Cheng, Li \& Qi, 2010) {\it Consider BCN (\ref{eq1}) with the algebraic form (\ref{eq4}). Assume the largest unobservable subspace $\mathcal{O}_c$ is a
regular subspace with $\{\,\tilde z_{\tilde s+1},\,\tilde z_{\tilde
s+2},\,\cdots,\,\tilde z_n\,\}$ as its basis. Then, under the coordinate transformation $\tilde z=\tilde T x$, the BCN (\ref{eq4}) becomes
\begin{shrinkeq}{-2ex}
\begin{equation}
\label{eq6}
\begin{array}{rcl}
\tilde z^{[1]}(t+1)& =& \tilde G_1u(t)\tilde z^{[1]}(t), \\
\tilde z^{[2]}(t+1)& =& \tilde G_2u(t)\tilde z(t),\\
y(t)& =& \tilde M\tilde z^{[1]}(t),
\end{array}
\end{equation}
\end{shrinkeq}
where $\tilde z^{[1]}=\ltimes_{i=1}^{\tilde s}\tilde z_i$,
$\tilde z^{[2]}=\ltimes_{i={\tilde s}+1}^n{\tilde z}_i$, $\tilde
G_1\in \mathcal{L}_{2^{\tilde s}\times 2^{\tilde s+m}}$, $\tilde
G_2\in \mathcal{L}_{2^{n-{\tilde s}}\times 2^{n+m}}$ and $\tilde
M\in \mathcal{L}_{2^{p}\times 2^{\tilde s}}$. The decomposition form (\ref{eq6}) is called the {\it observable normal form}
of (\ref{eq4}).}
\par
A natural question is whether the defined maximum decomposition
with respect to outputs and the observable normal form proposed by
Cheng, Li and Qi (2010) are the same one. Actually, with the
regularity assumption, they are the same one.
\par
For the basic concepts of the state-space method mentioned above, please refer to Cheng \& Qi (2010) and Cheng, Li, \& Qi (2010b). Comparing (\ref{eq5}) with (\ref{eq6}), we find that the maximum decomposition with respect to outputs and the observable normal form have the same structure.
\par
{\bf Proposition 3.} {\it Assume that the largest
unobservable subspace $\mathcal{O}_c$ is regular, then (\ref{eq5}) is
a maximum decomposition with respect to outputs if and only if it is an observable
normal form described by Proposition 2.}\\
{\bf Proof.} Assume that (\ref{eq5}) is a maximum decomposition with respect to outputs of order $n-s$. We first prove that $s=\tilde s$. On the one hand, by the definition of the maximum decomposition with respect to outputs,
we have $n-s\geq n-\tilde s$, i.e. $s\leq \tilde s$. On the other hand, from the definition of largest unobservable subspace, it follows that
\begin{shrinkeq}{-1ex}
\begin{equation}
\label{eqCc}
z_{s+1}, \cdots, z_n\in
\mathcal{O}_c=\mathcal{F}(\tilde z_{\tilde s+1}, \tilde z_{\tilde
s+2}, \cdots, \tilde z_n).
\end{equation}
\end{shrinkeq}
Thus, (\ref{eqCc}) implies that $n-s\leq n-\tilde s$, namely $s\geq \tilde s$. Therefore, we have proved that $s=\tilde s$, which implies that (\ref{eq6}) is a maximum decomposition with respect to outputs of  system (\ref{eq4}). Conversely, by $\tilde s=s$, (\ref{eqCc}) and Theorem 13 of Cheng, Li, \& Qi (2010), we obtain that $\{z_{ s+1}, \cdots, z_n\}$ is a regular basis of $\mathcal{O}_c$, that is, (\ref{eq5}) is an observable normal form. \hspace{\fill} \qed

\par
{\bf Remark 1}. {\it
Proposition 3 implies that the decomposition with respect to outputs is a generalization of the observable norm form. From the proof of Proposition 3, we see that, if the largest unobservable subspace $\mathcal{O}_c$ is regular, then $z^{[2]}$ and $\tilde z^{[2]}$ can be logically expressed in terms of each other,
i.e. there exists a logical matrix $R$ such that $(\textbf{1}_{2^{s}}^{\mathrm{T}}\otimes
I_{2^{n-s}})T=R\tilde T^{\mathrm{T}}(\textbf{1}_{2^{s}}\otimes
I_{2^{n-s}})$, i.e.
\begin{equation}
\label{eqR}
R=\frac{1}{2^s}(\textbf{1}_{2^{s}}^{\mathrm{T}}\otimes
I_{2^{n-s}})T\tilde T^{\mathrm{T}}(\textbf{1}_{2^{s}}\otimes
I_{2^{n-s}}).
\end{equation}
If there exist two maximum decompositions with respect
to outputs described by (\ref{eq5}) and (\ref{eq5}) respectively such that the right side of $(\ref{eqR})$
is not a logical matrix, then the largest uncontrollable subspace
is not regular.}
\section {Algebraic conditions for the decomposability with respect to outputs}

In this section, based on the definition of the decomposability with respect to outputs, we derive some equivalent algebraic conditions.

\par

{\bf Lemma 1.}{\it Assume that $M_1,\ M_2,\ \cdots, M_l\in
\mathbf{R}_{m\times n}$ are non-negative matrices satisfying
$\mathbf{1}_{m}^\mathrm{T}M_k=\mathbf{1}_{n}^\mathrm{T}$ for every $k=1,2,\cdots,l$. If $M_1+ M_2+\cdots+M_l=lG$, with $G$ being a logical matrix, then $M_1=M_2=\cdots =M_l=G$.}

\par

Swap matrix $W_{[m,n]}$ is an $mn\times mn$ logical matrix, defined as $W_{[m,n]}=[I_n\otimes \delta_m^1,I_n\otimes \delta_m^2,\cdots,I_n\otimes \delta_m^m].$
\par
{\bf Lemma 2.} (Cheng, Qi, \& Li, 2011) {\it Let $W_{[m,n]}\in \mathbf{R}_{mn\times mn}$ be a swap matrix. Then $W_{[m,n]}^\mathrm{T}=W_{[m,n]}^\mathrm{-1}=W_{[n,m]}$ and $W_{[m,1]}=W_{[1,m]}=I_m$, where $I_m$ is an identity matrix.}
\par
{\bf Lemma 3.} (Cheng, Qi, \& Li, 2011) {\it Let $A\in \mathbf{R}_{m\times n}$, $B\in \mathbf{R}_{p\times q}$. Then $W_{[m,p]}(A\otimes B)W_{[q,n]}=(B\otimes A)$.}

\par

{\bf Theorem 1.} {\it Consider BCN (\ref{eq1}) with the algebraic form (\ref{eq4}). Let $L=[L_1, L_2, \cdots, L_{2^m}]$, where $L_i\in \mathcal{L}_{2^n\times 2^{n}}$. Then the following statements are equivalent:\\
$(\mathrm{i})$ the system (\ref{eq1}) is decomposable with respect to outputs with order $n-s$; \\
$(\mathrm{ii})$ there exist a permutation matrix $T\in \mathcal{L}_{2^n\times 2^{n}}$, logical matrices $G_1\in \mathcal{L}_{2^s\times 2^{m+s}}$ and $M\in \mathcal{L}_{2^p\times 2^{s}}$ such that
\vspace{-0.3cm}
\begin{eqnarray}
\label{eq7}
  (I_{2^s}\otimes \mathbf{1}_{2^{n-s}}^{\mathrm{T}})TL(I_{2^m}\otimes T^{\mathrm{T}}) &=& G_1(I_{2^{m+s}}\otimes \mathbf{1}_{2^{n-s}}^{\mathrm{T}}), \\
\label{eq7-0}
  HT^{\mathrm{T}} &=& M(I_{2^{s}}\otimes \mathbf{1}_{2^{n-s}}^{\mathrm{T}}).
\end{eqnarray}
$(\mathrm{iii})$ there exist a permutation matrix
$T\in \mathcal{L}_{2^n\times 2^{n}}$,
$G_1\!=\![G_{11},G_{12}, \cdots, G_{12^m}]\in \mathcal{L}_{2^s\times 2^{m+s}}$
and $M\in \!\mathcal{L}_{2^p\times 2^{s}}$ such that
\begin{equation}
\label{eq8}
 QL_i=G_{1i}Q \ \ \mbox{ and }
  H=MQ
\end{equation}
hold for each $i=1,2,\cdots,2^m$, where $Q=(I_{2^s}\otimes \mathbf{1}_{2^{n-s}}^{\mathrm{T}})T$ and $G_{1i}\in \mathcal{L}_{2^s\times 2^{s}}$;\\
$(\mathrm{iv})$ there exists a permutation matrix $T\in \mathcal{L}_{2^n\times 2^{n}}$ such that
\begin{equation}
\label{eq9}
QL_iQ^\mathrm{T}/2^{n-s} \ \ \mbox{ and }
HQ^\mathrm{T}/2^{n-s}
\end{equation}
are logical matrices, where $Q=(I_{2^s}\otimes \mathbf{1}_{2^{n-s}}^{\mathrm{T}})T$.} \\
{\bf Proof}. (i) $\Leftrightarrow$ (ii) Assume that BCN (\ref{eq1})
is decomposable with respect to outputs of order $n-s$. Then
(\ref{eq3}) can be converted in (\ref{eq5}). By (\ref{eq3}), we have
\begin{eqnarray}
\label{eq10}
z^{[1]}(t+1)&=&(I_{2^{s}}\otimes\textbf{1}_{2^{n-s}}^{\mathrm{T}} )z(t+1) \nonumber\\
&=&(I_{2^{s}}\otimes\textbf{1}_{2^{n-s}}^{\mathrm{T}} )TL(I_{2^{m}}\otimes T^{\mathrm{T}})u(t)z(t),\\
\label{eq10-0}
y(t)&=&HT^{\mathrm{T}}z^{[1]}(t)z^{[2]}(t).
\end{eqnarray}
By (\ref{eq5}), we have
\begin{eqnarray}
\label{eq10-1}
z^{[1]}(t+1)&=&G_1(I_{2^{m+s}}\otimes \mathbf{1}_{2^{n-s}}^{\mathrm{T}})u(t)z(t),\\
\label{eq10-2}
y(t)&=&M(I_{2^{s}}\otimes \mathbf{1}_{2^{n-s}}^{\mathrm{T}})z^{[1]}(t)z^{[2]}(t).
\end{eqnarray}
From (\ref{eq10}) and (\ref{eq10-1}), we get (\ref{eq7}). From (\ref{eq10-0}) and (\ref{eq10-2}), we get (\ref{eq7-0}). Conversely, with the similar procedure, it follows from (\ref{eq7}) and (\ref{eq7-0}) that (\ref{eq5}) holds. Thus (i) is proved.\\
(ii) $\Rightarrow$ (iii) Multiplying (\ref{eq7}) on
the right by $I_{2^{m}}\otimes T$ yields
\begin{eqnarray}
QL&=&G_1(I_{2^m}\otimes((I_{2^{s}}\otimes \mathbf{1}_{2^{n-s}}^{\mathrm{T}})T))\nonumber \\
&=& [G_{11},\ G_{12},\ \cdots, G_{12^m}] (I_{2^m}\otimes Q)\nonumber \\
&=&[G_{11}Q,\ G_{12}Q, \cdots, G_{12^m}Q],
\end{eqnarray}
which implies $QL_i=G_{1i}Q$. Multiplying (\ref{eq7-0}) on
the right by $T$ gives $H=MQ$.\\
(iii) $\Rightarrow$ (iv) A straightforward calculation shows that
$$QQ^{\mathrm{T}}=(I_{2^s}\otimes \mathbf{1}_{2^{n-s}}^{\mathrm{T}})TT^{\mathrm{T}}(I_{2^{s}}\otimes\mathbf{1}_{2^{n-s}})
=2^{n-s}I_{2^{s}}.$$
Thus, it follows from (\ref{eq8}) that
\begin{equation}
\label{eq9-1}
G_{1i}=QL_iQ^\mathrm{T}/2^{n-s} \ \ \mbox{ and }
M=HQ^\mathrm{T}/2^{n-s}.
\end{equation}
Thus the matrices in (\ref{eq9}) are logical matrices.\par
(iv) $\Rightarrow$ (iii) Denote the logical matrices in (\ref{eq9}) by $G_{1i}$ and $M$ respectively. Let
\begin{shrinkeq}{-1ex}
\begin{equation}
\label{eq11}
QL_iT^{\mathrm{T}}W_{[2^{n-s},2^s]}=[P_1, P_2, \cdots, P_{2^{n-s}}],
\end{equation}
\end{shrinkeq}
where $P_i\in \mathcal{L}_{2^{s}\times 2^{s}}$ are non-negative matrices.
By (\ref{eq11}), Lemmas 2 and 3, we have
\begin{equation}
\label{eq12}
\begin{array}{rcl}
 G_{1i}&=&QL_iQ^\mathrm{T}/2^{n-s} \nonumber \\
  &=&[P_1,\!P_2,\!\cdots,\!P_{2^{n-s}}]W_{[2^s,2^{n-s}]}TQ^\mathrm{T}/2^{n-s}\nonumber \\
  &=&[P_1, P_2, \cdots, P_{2^{n-s}}](\textbf{1}_{2^{n-s}}\otimes I_{2^{s}})/{2^{n-s}},
\end{array}
\end{equation}
that is,
\begin{equation}
\label{eq13}
\sum_{k=1}^{2^{n-s}}P_k=2^{n-s}G_{1i}.
\end{equation}
Multiplying (\ref{eq11}) on the left by $\mathbf{1}_{2^{s}}^{\mathrm{T}}$ yields
\begin{equation}
\label{eq13-1}
\mathbf{1}_{2^{s}}^{\mathrm{T}}[P_1, P_2, \cdots, P_{2^{n-s}}]=\mathbf{1}_{2^{n}}^{\mathrm{T}} W_{[2^{n-s},2^s]}=\mathbf{1}_{2^{n}}^{\mathrm{T}},
\end{equation}
namely, $\mathbf{1}_{2^{s}}^{\mathrm{T}}P_k=\mathbf{1}_{2^{s}}^{\mathrm{T}}$ for every $k=1,2,\cdots, 2^{n-s}$.
Thus, from Lemma 1 and (\ref{eq13}), it follows that $P_k=G_{1i}$.
Considering (\ref{eq11}), we have
\begin{equation}
\label{eq13-2} QL_iT^{\mathrm{T}}W_{[2^{n-s},2^s]}=G_{1i}(\textbf{1}^{\mathrm{T}}_{2^{n-s}}\otimes I_{2^{s}}).
\end{equation}
Thus
\begin{equation}
\label{eq13-3} QL_iT^{\mathrm{T}}=G_{1i}( I_{2^{s}}\otimes\textbf{1}^{\mathrm{T}}_{2^{n-s}}),
\end{equation}
which implies $QL_i=G_{1i}Q$.
With a same procedure above, we obtain $H=MQ$.\\
(iii) $\Rightarrow$ (ii) The proof is trivial.
\hspace{\fill} \qed
\par

\section{A graphical condition for the decomposability with respect to outputs }
Consider BCN (\ref{eq4}) with $L\!=\![L_1,L_2,\cdots,\!L_{2^m}\!]\!\in \!\mathcal{L}_{2^n\!\times\! 2^{n+m}}$. Each logical matrix $L_j$ can be regarded as an adjacency matrix of a directed
graph $\mathcal{G}_j$. The vertex set of $\mathcal{G}_j$ is $A=\{1,2,3,\cdots,2^n\}$, and $\mathcal{G}_j$ has a directed edge $(q,k)$ if and only if
$(L_j)_{kq}\ne 0$. We say that $k$ is an {\it out-neighbor} of $q$ with respect to $\mathcal{G}_j$
if $(L_j)_{kq}\ne 0$. We denote the {\it out-neighborhood} of set $S$ in graph $\mathcal{G}_j$ by $\mathcal{N}^j(S)$.

{\bf Definition 1.} Let $A$ be the vertex set of a graph $\mathcal{G}$, and $\Phi_l,l=1,2,\cdots,\mu$ be subsets of $A$.  $\{\Phi_l\}_{l=1}^\mu$ is
called a {\it vertex partition} of $A$, if $\cup_{l=1}^\mu\Phi_l=A$ and $\Phi_i\cap
\Phi_j=\emptyset$ for any $i\ne j$. A vertex
partition $\{S_l\}_{l=1}^\mu$ of $A$ is called an {\it
equal vertex partition} if
$|S_l|=|A|/\mu$ for every $l=1,2,\cdots, \mu$.

\par


In Zou \& Zhu (2014), the concept of {\it perfect equal vertex partition} (PEVP) is proposed for the decomposition with respect to inputs. Here, we further extend this concept to deal with the decomposition with respect to outputs.
\par
{\bf Definition 2}. {\it Consider the set $\{\mathcal{G}_j\}_{j=1}^{2^m}$ composed of $2^m$ digraphs with the common vertex set $A$. Assume that   different colors are assigned to all the
vertices of $A$. An equal vertex partition $\mathcal{S}=\{S_l\}_{l=1}^\mu$ of $A$ is called a {\it common concolorous perfect equal vertex partition} (CC-PEVP) of $\{\mathcal{G}_j\}_{j=1}^{2^m}$ if
\\
$(\mathrm{i})$ $\mathcal{S}$ is perfect for each $j=1,2,\cdots,2^m$, i.e. for any given $l$ and $j$, there exists an $\alpha_l^j$ such that $\mathcal{N}^j(S_l)\subset S_{\alpha_l^j}$ ;\\
$(\mathrm{ii})$ for any given $l=1,2,\cdots,\mu$, all the vertices in $S_l$ have the same color.}
\par

\begin{figure}
\begin{center}
\includegraphics[width=5cm]{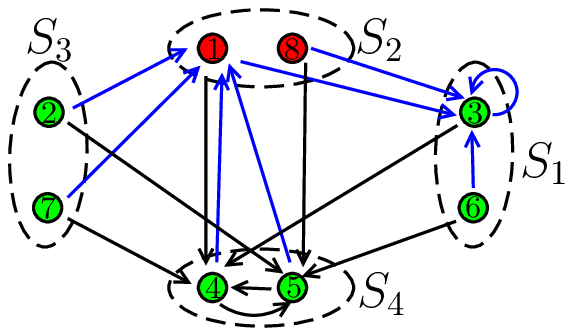}    
\\ Fig.1. CC-PEVP of two graphs $\mathcal{G}_1$ and $\mathcal{G}_2$.
\label{fig1}                                 
\end{center}                                 
\end{figure}

We give a simple example shown in Fig. 1 to explain Definition 2 intuitively. In Fig. 1, each vertex has red or green color. The digraph with blue edges is denoted by $\mathcal{G}_1$ and the other one with black edges is denoted by $\mathcal{G}_2$. Obviously, $\{S_1,S_2,S_3,S_4\}$ shown in Fig.1 forms an equal vertex partition and each $S_l$ contains vertices of the same color. Moreover, from Fig.1, we have
\vspace{-2mm}
$$
\begin{array}{c}
\mathcal{N}^1\!(\!S_1\!)\!=\!\mathcal{N}^1\!(\!S_2\!)\!=\!\{3\}\!\subset \!S_{1},\ \ \mathcal{N}^1\!(\!S_3\!)\!=\!\mathcal{N}^1\!(\!S_4\!)\!=\!\{1\}\!\subset\!S_{2},\\ \mathcal{N}^2\!(\!S_1\!)\!=\!\mathcal{N}^2\!(\!S_2\!)\!=\! \mathcal{N}^2\!(\!S_3\!)\!=\!\mathcal{N}^2\!(\!S_4\!)\!=\!\{4,\ 5\}\!\subset\!S_{4}.\
\end{array}
$$
Therefore, by Definition 2, $\{S_1,S_2,S_3,S_4\}$ is a CC-PEVP of $\{\mathcal{G}_1,\ \mathcal{G}_2\}$.
\par
Based on Definition 2, we propose an equivalent graphical condition for the decomposability with respect to outputs in the following theorem.
\par
{\bf Theorem 2.} {\it Consider the BCN (\ref{eq1}) with algebraic form (\ref{eq4}). Let $L=[L_1,L_2,\cdots,L_{2^m}]$ with each $L_j\in \mathcal{L}_{2^n\times 2^n}$. Denote by $A$ the common vertex set of $\{\mathcal{G}_j\}_{j=1}^{2^m}$, where each $\mathcal{G}_j$ is induced by $L_j$. Color all the vertices in $A$ in such a way that any two vertices $\mu$ and $\lambda$ are of the same color if and only if $H\delta_{2^n}^\mu=H\delta_{2^n}^\lambda$. Then
BCN (\ref{eq1}) is  decomposable with respect to outputs with
order $n-s$ if and only if $\{\mathcal{G}_j\}_{j=1}^{2^m}$ has a CC-PEVP  $\{S_l\}_{l=1}^{2^{s}}$ with $|S_l|=2^{n-s}$.}
\par
Before the proof of Theorem 2, we give an intuitive explanation on the motivation.
From the decomposition form (\ref{eq3}), it follows that, as  $z_{1},\cdots,z_s$ are fixed,
the set
\begin{equation}
\label{set} S_{z_{1}\cdots
z_s}\!\!=\!\!\{\!(z_1,\!\cdots,\!z_s,\!z_{s+1},\!\cdots,\!z_n)\ | \ z_j\!\in\!\{0,\!1\},s+1\!\leq \!j\!\leq\!n\!\}\nonumber
\end{equation}
has $2^{n-s}$ states and these states have the same output. Then the family
\begin{shrinkeq}{-1ex}
\begin{equation}
\label{family}\{S_{z_{1}\cdots z_s}\ |\ z_j\in \{0,\ 1\}, \
j=1,2,\cdots,s\},
\end{equation}
\end{shrinkeq}
forms a concolorous equal partition of all the $2^n$ states. From the algebraic representation (\ref{eq5}), it follows that  $z^{[1]}(t+1) = G_{1j}z^{[1]}(t)$ for any fixed $u(t)=\delta_{2^m}^j$, i.e. the state $z(t)$ transmits from one $S_{z_{1}\cdots z_s}$ to another for any given graph $\mathcal{G}_j$. This is just why we propose the concept of CC-PEVP and try to reveal the relationship between CC-PEVP and the decomposition with respect to outputs.\par
\par
{\bf Proof of Theorem 2}. (Necessity) By Theorem 1, there exist a
permutation matrix $T\in \mathcal{L}_{2^n\times 2^{n}}$ and
logical matrices $G_{1j}\in \mathcal{L}_{2^{s}\times 2^{s}}(j=1,2,\cdots,2^m)$, $M\in \mathcal{L}_{2^{p}\times 2^{s}}$ such
that the equalities in (\ref{eq8}) hold for each $j=1,2,\cdots,2^m$. Set
\begin{equation}
\label{eq15}
Q=(I_{2^{s}}\otimes \textbf{1}_{2^{n-s}}^{\mathrm{T}} )T=\delta_{2^{s}}[i_1,\cdots,i_{2^n}].
\end{equation}
Thus, (\ref{eq8}) can be rewritten as
\begin{eqnarray}
\label{eq17-3}
\hspace{1cm} \delta_{2^{s}}[i_1,\cdots,i_{2^n}]L_j&=&G_{1j}\delta_{2^{s}}
[i_1,\cdots,i_{2^n}],\\
\label{eq18}
\hspace{1cm}  H&=&M\delta_{2^{s}}[i_1,\cdots,i_{2^n}].
\end{eqnarray}
Let $S_l=\{q| i_q=l\}$. Then $\{S_l\}_{l=1}^{2^{s}}$ is an equal vertex
partition of $A$ with $|S_l|=2^{n-s}$. For any $l$, we have that
\begin{equation}\label{eq20}
  \exists \ \alpha^j_l,\ \beta_l, \mbox{ s.t. }  G_{1j}\delta_{2^{s}}^l=\delta_{2^{s}}^{\alpha^j_l},\ M\delta_{2^{s}}^l=\delta_{2^{p}}^{\beta_l}.
\end{equation}
For any $k\in
\mathcal{N}^j(S_l)$, there exists $q\in S_l$, i.e. $i_q=l$, such that
$\textrm{Col}_q(L_j)=\delta_{2^n}^k$. Then, by (\ref{eq20}), we have
\begin{equation}
\label{eq21}
\delta_{2^{s}}^{i_k}\!=\!\delta_{2^{s}}\![i_1,\!\cdots,\!i_{2^n}\!]\textrm{Col}_q(L\!_j\!)
\!=\!G_{\!1\!j}\delta_{2^{s}}^{i_{q}}\!=\!
G_{\!1\!j}\delta_{2^{s}}^{l}\!\!=\!
\delta_{2^{s}}^{\alpha^j_l}.
\end{equation}
Thus $i_k=\alpha^j_l$, which implies $k\in S_{\alpha^j_l}$. Therefore, we have $\mathcal{N}^j(S_l)\subset S_{\alpha^j_l}$. Moreover, for every $k\in S_l$ ($i_k=l$), $$H\delta_{2^n}^k=M\delta_{2^{s}}^{i_k}=M\delta_{2^{s}}^{l}=
\delta_{2^{p}}^{\beta_l},$$ which implies that all the vertices in $S_l$ have the same color.\par
(Sufficiency) Since $\{S_l\}_{l=1}^{2^{s}}$ is a CC-PEVP, by Definition 1 and the definition of the coloring, we have that
\begin{equation}
\label{eqeq29}
\forall \ 1\leq l\leq 2^s,\ 1\leq j\leq 2^m,\ \exists \ \alpha_l^j,\ \mbox{s.t.}\  \mathcal{N}^j(S_l)\subset S_{\alpha_l^j}
\end{equation}
and
\begin{equation}
\label{eqH}
\forall \ 1\leq l\leq 2^s,\  \exists \ \beta_l,\ \mbox{s.t.}\  \forall \ q\in S_l,\ H\delta_{2^n}^q=\delta_{2^p}^{\beta_l}.
\end{equation}
For any $q\in S_l$ ($l=1,2,\cdots,2^s$), let $i_q=l$. We denote $Q=\delta_{2^{s}}[i_1,\cdots,i_{2^n}]$. Since $|S_l|=2^{n-s}$, there exists a permutation matrix $T$ such that
$(I_{2^{s}}\otimes \textbf{1}_{2^{n-s}}^{\mathrm{T}})T=Q$. For any $l\in \{1,2,\cdots,2^s\}$, we have
\vspace{-0.3cm}
\begin{eqnarray}
\label{eq25}
\mbox{Col}_l(\!QL_jQ^\mathrm{T}\!)&=&QL_j \mbox{Col}_l(\!Q^\mathrm{T}\!)\!=\!QL\!_j\!\!\sum_{q\in S_l}\!\delta_{2^{n}}^{q}\nonumber \\
&=&\sum_{q\in S_l}\!\!Q\mbox{Col}_q(L\!_j).
\end{eqnarray}
Since $L\!_j$ is a logical matrix, we can let $\mbox{Col}_q(L\!_j)=\delta_{2^n}^{k_q}$, i.e. $k_q\in \mathcal{N}^j(\{q\})$. Thus it follows from (\ref{eq25}) and (\ref{eqeq29}) that
\begin{eqnarray}
\label{eq25-5}
\mbox{Col}_l(\!QL_jQ^\mathrm{T}\!)\!&=&\!\!\sum_{q\in S_l}Q\delta_{2^n}^{k_q}=\sum_{q\in S_l}\delta_{2^{s}}^{i_{k_q}}\nonumber \\
\!&=&\sum_{q\in S_l}\delta_{2^{s}}^{\alpha^j_l}\ \  (k_q\in \!N^j(\{q\})\!\subset \!N^j(S_l)\!\subset \!S_{\alpha^j_l})\nonumber \\
&=&2^{n-s}\delta_{2^{s}}^{\alpha^j_l}.
\end{eqnarray}
From (\ref{eqH}), it follows that
\begin{eqnarray}
\label{eq25-6}
\mbox{Col}_l(HQ^\mathrm{T})&=&H\mbox{Col}_l(Q^\mathrm{T})=\sum_{q\in S_l}H\delta_{2^{n}}^{q} \nonumber\\
& =&\sum_{q\in S_l}\delta_{2^{s}}^{\beta_l}=2^{n-s}\delta_{2^{s}}^{\beta_l}.
\end{eqnarray}
By (\ref{eq25-5}) and (\ref{eq25-6}), both $QL_jQ^\mathrm{T}/2^{n-s}$ and $HQ^\mathrm{T}/2^{n-s}$ are
logical matrices. Therefore, by (iv) of Theorem 1, the sufficiency is
proved.\hspace{\fill} \qed

\par

From the sufficiency proof of Theorem 2, we see that, if a
CC-PEVP $\{S_l\}_{l=1}^{\!2\!^s}$ is given, a constructive procedure
to calculate the transformation matrix $T$ is obtained as follow:
\\
(i) let $i_q=l$ for any $l=1,2,\cdots,2^s$ and $q\in S_l$; \\
(ii) let $Q=\delta_{2^{s}}[i_1,\cdots,i_{2^n}]$;\\
(iii) compute permutation matrix $T$ from
$$(I_{2^{s}}\otimes \textbf{1}_{2^{n-s}}^{\mathrm{T}})T=Q.$$
To display the effectiveness of Theorem 2, we
reconsider Example 10.3 of Cheng, Qi, \& Li (2011) and construct
the coordinate transformation using the graphical method.\par
{\bf Example 1.} Consider the following system
\begin{shrinkeq}{-2ex}
\begin{equation}
\label{eq25-00}
    \begin{array}{ll}
\!\!    x_{1}(t\!+\!1)\!\!=\!\!x_{3}(t)\vee u(t), \\
\!\!    x_{2}(t\!+\!1)\!\!=\!\!(x_{1}(t)\!\!\wedge\!\neg
                    x_3(t))\!\!\vee\!\!(\neg x_{1}(t)\!\!\wedge\!\!(x_{3}(t)\!\!\leftrightarrow\!\! u(t))),\\
\!\!    x_{3}(t\!+\!1)\!\!=\!\!x_{3}(t)\rightarrow u(t),
\\
  \quad\!\! y(t)=(x_1(t)\leftrightarrow x_3(t))\rightarrow(x_2(t)\bar{\vee} x_3(t)).
   \end{array}
\end{equation}
\end{shrinkeq}
Let $x(t)=x_{1}(t)x_{2}(t)x_{3}(t)$. Then we have
$$x(t+1)=Lu(t)x(t), \ \ y(t)=Hx(t),$$
where $L=[L_{1}, L_{2}]\in \mathcal{L}_{8\times 16}$, with $$L_{1}=\delta_{8}[3~1~3~1~1~3~1~3],\ \ L_{2} =\delta_{8}[4~5~4~5~4~5~4~5]$$ and $$H =\delta_{2}[2~1~1~1~1~1~1~2].$$
The digraphs $\mathcal{G}_1$ and $\mathcal{G}_2$ corresponding to BCN (\ref{eq25-00}) is just shown in Fig.1. The vertex partition $S_i\ (i=1,2,3,4)$ given in Fig.1 is a CC-PEVP of $\{\mathcal{G}_1,\ \mathcal{G}_2\}$. Thus BCN (\ref{eq25-00}) is decomposable with respect to outputs of order $1$. Let
\begin{equation*}
(I_{4}\otimes\textbf{1}_{2}^{\mathrm{T}})T=Q=\delta_{4}[2~3~1~4~4~1~3~2].
\end{equation*}
It follows that $T=\delta_{8}[3~6~1~8~7~2~5~4]$. The coordinate
transformation matrix $T$ is the same as that given in Example
10.3 of Cheng, Qi, \& Li (2011).
Therefore, the decomposition with respect to outputs is obtained as
\begin{equation*}
\begin{array}{l}
  \left\{
   \begin{array}{ll}
\!\!    z_{1}(t+1)=u(t), \\
\!\!    z_{2}(t+1)=z_{1}(t)\wedge u(t), \\
\!\!    z_{3}(t+1)=z_{3}(t)\rightarrow u(t).
   \end{array}
 \right. \\
  \quad\!\! y(t)=z_1(t)\rightarrow z_2(t).
\end{array}
\end{equation*}
\section {Searching a CC-PEVP}
To achieve the decomposition with respect to outputs of the
maximum order $n-s$, we need to find a CC-PEVP
$\{S_l\}_{l=1}^{2^s}$ of the minimum $s$. Since the number of
all the equal vertex partitions is finite, a straightforward
method is to check whether each one is a CC-PEVP or not, but
the computational burden of this method is very heavy. To reduce the
computational burden, we investigate some necessary conditions for the
existence of a CC-PEVP.
\par
For BCN (\ref{eq1}) with the algebraic form (\ref{eq4}), let
\begin{equation}
\label{Ri} R^{j_1j_2\cdots j_r}_k=\{q|\
\mbox{Col}_q(HL_{j_1}L_{j_2}\cdots L_{j_r})=\delta_{2^p}^k\}
\end{equation}
for any $k\!=\!1,2,\cdots,2 ^{p}$, $r\!\geq \!0$ and $1\!\leq\!
j_1,j_2,\cdots, j_r\!\leq \!2^m$. Here, $r=0$ means that all
the $j_i$ vanish, i.e.
$$R_k=\{q|\ \mbox{Col}_q(H)=\delta_{2^p}^k \}.$$
From (\ref{Ri}), we get the following result.
\par
{\bf Lemma 4}. {\it For any given $j_1$, $j_2$, $\cdots$ $j_r$, the
family $\mathcal{R}^{j_1j_2\cdots j_r}=\{R^{j_1j_2\cdots
j_r}_k\}_{k=1}^{2^p}$ is a partition of $A=\{1,2,\cdots,2^n\}$.
For any $q_1,\ q_2\in A$, they are in the same
$R^{j_1j_2\cdots j_r}_k$ if and only if
$$\mbox{Col}_{q_1}(HL\!_{j_1}\!\cdots
L\!_{j_r})=\mbox{Col}_{q_2}(HL\!_{j_1}\!\cdots L\!_{j_r}).$$}
\par
Given a partition $\mathcal{R}$ of set $A$, for any $q_1, q_2\in A$, we say that $q_1$ is equivalent to $q_2$, denoted by $q_1 \overset{\mathcal{R}}{\sim} q_2$, if and only if there exists a $R\in \mathcal{R}$ such that $q_1, q_2\in R$. Conversely, each equivalence relation on $A$ induces an associated partition $\mathcal{R}$ of set $A$. Therefore, by Lemma 4, we have
\begin{equation}
\label{Eqiv}q_1 \!\!\overset{\mathcal{R}^{j\!_1\cdots j\!_r}}{\sim} \!\!q_2 \Leftrightarrow \mbox{Col}_{q_1}\!(\!H\!L\!_{j_1}\!\!\cdots
\!L\!_{j_r}\!)\!=\!\mbox{Col}_{q_2}\!(\!H\!L\!_{j_1}\!\!\cdots\!
L\!_{j_r}\!).
\end{equation}
\ \par
Let us briefly recall some basic and well-known terms, notions and facts, without proofs, concerning partitions and equivalence relations. For details, please refer to Bor\r{u}vka (1974) and Pot\r{u}{\v c}ek (2014).\par
{\bf Definition 3}. {\it Let $\mathcal{U}$ and $\mathcal{V}$ be partitions of a set $A$. Assume that, for every $U\in \mathcal{U}$, there exists $V\in \mathcal{V}$ such that $U \subset V$. Then partition $\mathcal{U}$ is said to be a {\it refinement} of $\mathcal{V}$, which is denoted by $\mathcal{U}\sqsubset\mathcal{V}$.}
\par
{\bf Definition 4}. {\it A {\it meet} of two partitions $\mathcal{U}$ and $\mathcal{V}$ of a set $A$, denoted by $\mathcal{U}\sqcap\mathcal{V}$,
is a set of all intersections $U\cap V$, where $U \in \mathcal{U}$ and  $V \in \mathcal{V}$.}
\par
{\bf Lemma 5}. {\it Let $\mathcal{V}_1$ and $\mathcal{V}_2$ be two
partitions of a set $A$ and let
$\mathcal{V}=(\mathcal{V}_1\sqcap\mathcal{V}_2)\setminus
\{\emptyset\}$. Then $\mathcal{V}$ is a partition of $A$
satisfying $\mathcal{V}\sqsubset \mathcal{V}_1$ and
$\mathcal{V}\sqsubset \mathcal{V}_2$. Moreover, if $\mathcal{U}$
is a partition of $A$ such that $\mathcal{U}\sqsubset
\mathcal{V}_1$ and $\mathcal{U}\sqsubset \mathcal{V}_2$, then
$\mathcal{U}\sqsubset \mathcal{V}$.
}

{\bf Definition 5}. {\it The partition $\mathcal{V}$ given in Lemma 5 is called the {\it greatest common
refinement} of the partitions $\mathcal{V}_1$ and $\mathcal{V}_2$ denoted by $\mathbf{gcr}(\mathcal{V}_1,\ \mathcal{V}_2)$.}
\par
{\bf Definition 6}. {\it The {\it greatest common
refinement} of the partitions $\mathcal{V}_1$, $\mathcal{V}_2$, $\cdots$,
$\mathcal{V}_\mu$ of $A$ is inductively defined as \\
$\mathbf{gcr}(\mathcal{V}_1,\ \mathcal{V}_2,\cdots, \mathcal{V}_{\mu})\!=\!\mathbf{gcr}(\mathbf{gcr}(\mathcal{V}_1,\ \mathcal{V}_2,\cdots, \mathcal{V}_{\mu-1}), \mathcal{V}_\mu$).}


{\bf Lemma 6.} {\it Let $\mathcal{W}=\mathbf{gcr}(\mathcal{V}_1,\
\mathcal{V}_2,\cdots, \mathcal{V}_{\mu})$. Then, for any $q_1,
q_2\in A$, we have
\begin{equation}
\label{Eqiv2}q_1 \!\overset{\mathcal{W}}{\sim} \!q_2
\Leftrightarrow (q_1 \!\overset{\mathcal{\mathcal{V}}_i}{\sim}
\!q_2,\ \ i=1,2,\cdots,\mu).
\end{equation}}
Now, let us come to our main result:\par
{\bf Proposition 4.} {\it Assume that BCN (\ref{eq4}) with
$L\!=\![L_1,\cdots,\!L_{2^m}\!]$ is decomposable with respect to
outputs with order $n-s$, where each $L_i\!\in
\!\mathcal{L}_{2^n\!\times\! 2^{n}}$. Denote by $\mathcal{G}_j$
the directed graph with the adjacency matrix $L_j$. Let
$\mathcal{S}\!\!=\!\!\{\!S_l\!\}_{l=1}^{2^{s}}\
(|S_l|\!\!=\!\!2^{n-s})$ be a CC-PEVP of
 $\{\!\mathcal{G}_j\!\}_{j=1}^{2^m}$. Then 
$$\mathcal{S}\sqsubset \mathcal{R}^{j_1\cdots j_r},\ \  \forall\  r\geq 0,\ \ \forall\ j_1,\cdots,j_r\in \ \{1,2,\cdots,
 2^m\}.$$}
\\ {\bf Proof.} For any $l=1,2,\cdots,2^s$ and $q\in S_l$, let $i_q=l$. Set  $Q=\delta_{2^{s}}[i_1,\cdots,i_{2^n}]$. Then by the sufficiency proof of Theorem 2, we obtain (iv) of
Theorem 1, which is equivalent to (iii) of Theorem 1. From (\ref{eq8}), we have $H\!=\!MQ$ and
\begin{equation}
\label{eqHLj} HL\!_{j_1}\!\cdots
L\!_{j_r}\!\!=\!MQL\!_{j_1}\!\cdots L\!_{j_r}
 \!\!=\!MG_{1\!j_1}\cdots G_{1\!j_r}\!Q
\end{equation}
for any $1\leq j_1,\cdots,j_r \leq 2^m$. Choose
$\beta^{j_1j_2\cdots j_r}_l$ such that
\begin{equation}
\delta^{\beta^{j_1j_2\cdots j_r}_l}_{2^p}=MG_{1\!j_1}\cdots
G_{1\!j_r}\delta^{l}_{2^s}.
\end{equation}
Now we are in a position to prove $S_l\subset R^{j_1j_2\cdots
j_r}_{\beta^{j_1j_2\cdots j_r}_l}.$ In fact, for every $k\in S_l$
($i_k=l$), it follows from (\ref{eqHLj}) that
\begin{eqnarray}
\label{ColHLj}
\mbox{Col}_k(HL\!_{j_1}\!\cdots L\!_{j_r})\!&=&\!MG_{1\!j_1}\cdots G_{1\!j_r}\!Q\delta_{2^n}^k\nonumber \\
&=&\!MG_{1\!j_1}\cdots G_{1\!j_r}\!\delta_{2^s}^{i_k}\nonumber\\
&=&\!MG_{1\!j_1}\cdots G_{1\!j_r}\!\delta_{2^s}^{l}\nonumber \\
 &=& \delta^{\beta^{j_1j_2\cdots j_r}_l}_{2^p}.
\end{eqnarray}
Thus it follows
from (\ref{Ri}) that $k\in R^{j_1j_2\cdots
j_r}_{\beta^{j_1j_2\cdots j_r}_l}$. \hspace{\fill} \qed
\par
From Proposition 4, it seems that an infinite number of partitions need to be
considered. But actually only a finite number of partitions are necessary due to the
following result.
\par
{\bf Lemma 7} (Cheng, Qi, \& Zhao, 2011). {\it Define a set of logical
matrices as
\begin{eqnarray}
\label{Omega}
&&\mathcal{H}_0\!=\!\{H\},\nonumber \\
&&\mathcal{H}_r\!=\!\{HL_{j_1}L_{j_2},\cdots,L_{j_r}|\
j_1,j_2,\cdots,j_r\!\in \!\{1,2,\cdots,2^m\}\},\nonumber
\end{eqnarray}
where $r=1,2,\cdots$. Then there exists the minimum $r^*$ such
that $\mathcal{H}_r \subset \cup_{r=0}^{r^*}\mathcal{H}_{r}$ for
any $r>r^*$.} \par
Denote
\begin{equation}
\label{eq26-1}
  \mathcal{H}=\mathcal{H}_0\cup\mathcal{H}_1\cup\cdots
  \cup\mathcal{H}_{r^*}
\end{equation}
and let $|\mathcal{H}|=\tau$.
For convenience, we denote the $\tau$ matrices in $\mathcal{H}$ by
$H_i=\delta_{2^p}[h_{i1},h_{i2},\cdots,h_{i2^n}]$ and the corresponding partitions by $\mathcal{R}_i$ ($i=1,2,\cdots,\tau$).
By Proposition 4 and Lemma 7, we
see that it is only needed to search a CC-PEVP
$\{S_l\}_{l=1}^{2^s}$ from the partitions $\mathcal{R}_i$ ($i=1,2,\cdots,\tau$). Let
\begin{equation}
\label{Omatrix}
\mathcal{O}=\left[\begin{array}{cccc}
                    h_{11}&h_{12}&\cdots&h_{12^n}\\
                    h_{21}&h_{22}&\cdots&h_{22^n}\\
                    \vdots &  \vdots &  \ddots &  \vdots \\
                    h_{\tau 1}&h_{\tau 2}&\cdots&h_{\tau 2^n}
                    \end{array}
                    \right],
\end{equation}
which is just the observability matrix proposed in Cheng \& Qi (2009).
\par
{\bf Corollary 1}. {\it Under the conditions of Proposition 4, we have
\begin{equation}
\label{Alg} \mathcal{S} \sqsubset \mathbf{gcr}(\mathcal{R}_1,\mathcal{R}_2,\cdots,\mathcal{R}_\tau)=:\mathcal{C}.
\end{equation}}
\ \\ By Corollary 1, we only need to search the CC-PEVP
$\mathcal{S}=\{S_l\}_{l=1}^{2^s}$ of the minimum $s$ from
partition $\mathcal{C}$.
\par
{\bf Proposition 5}. {\it For any $q_1,\ q_2\in A$, we have
\begin{equation}
\label{Ceqiv} q_1\overset{\mathcal{C}}{\sim} q_2 \ \Leftrightarrow
\ \mbox{Col}_{q_1}(\mathcal{O})=\mbox{Col}_{q_2}(\mathcal{O}),
\end{equation}
where partition $\mathcal{C}$ of the vertex set $A$ is shown in (\ref{Alg}).}\\
{\bf Proof}. By (\ref{Eqiv}), Lemma 6 and (\ref{Omatrix}), the
proposition is proved.
\par
{\bf Corollary 2.} {\it If there exists a $C\in \mathcal{C}$ with odd
cardinal, then BCN (\ref{eq4}) is undecomposable with respect to
outputs.}
\\
{\bf Proof}. By Corollary 1, $C$ is a union of some $S_l$. Since
$|S_l|=2^{n-s}$ for each $l$, we see that $2^{n-s}$ is a factor of
$|C|$. Considering that $|C|$ is odd, we obtain the order $n-s=0$.
Thus BCN (\ref{eq4}) is undecomposable with respect to outputs.
\par
{\bf Lemma 8} (Cheng \& Qi, 2009). {\it Assume that BCN (\ref{eq4}) is
globally controllable. Then it is observable if and only
if all the columns of $\mathcal{O}$ are distinct.}
\par
{\bf Corollary 3.} {\it Assume that BCN (\ref{eq4}) is globally
controllable. If (\ref{eq4}) is observable, then it is
undecomposable with respect to outputs.}
\\
{\bf Proof}. From Lemma 8, it follows that $|C|=1$ for each $C\in
\mathcal{C}$. Thus the corollary is proved by Corollary 2.
\par
{\bf Remark 2}. {\it The inverse of Corollary 3 is not correct. For
example, consider a Boolean network (\ref{eq4}) without control
($m=0$). Let $L=\delta_4[1, 2, 3, 1]$ and $H=\delta_2[1, 1, 1, 2]$.
Since $HL^i=\delta_2[1, 1, 1, 1]$ for any $i>0$, we have that
$$
\mathcal{O}=\left[\begin{array}{cccc}
                    1&1&1&2\\
                    1&1&1&1
                    \end{array}
                    \right],
$$
which implies that $\mathcal{C}=\{\{1,2,3\},\ \{4\}\}$. Thus the
BCN is unobservable and undecomposable. However, for the
traditional linear control systems, the observability and the
decomposability are equivalent. This is a property of
BCNs different from that of the traditional linear control systems. Hence, we prefer to call (\ref{eq3}) the decomposition with
respect to outputs instead of the observability decomposition.}
\par
To illustrate our method, we consider a practical example. \par
{\bf Example 2}. A shift register is a cascade of flip-flops:
\begin{eqnarray}
x_i(k+1)&=&x_{i+1}(k), \ \ (i=1,2,\cdots, n-1),\nonumber \\
 x_n(k+1)&=&u(k), \nonumber \\
 y(k)&=&x_1(k),
 \end{eqnarray}
where $x_i$ is the binary state of the $i$th flip-flop, $u$ the
input and $y$ the output. Multiplying all the equations yields
$$x(k+1)=Lx(k),\ \ \ y(k)=Hx(k),$$
where $L=[L_1, L_2]=\mathbf{1}_2^\mathrm{T}W_{[2,2^n]}$ with \vspace{-3mm}
\begin{eqnarray}
\hspace{1cm}L_1&=&\delta_{2^n}[1,3,5,\cdots,2^n\!\!-\!\!1,1,3,5,\cdots,2^n\!\!-\!\!1],\nonumber \\
\hspace{1cm}L_2&=&\delta_{2^n}[2,4,6,\cdots,2^n,\ \ \  2,4,6,\cdots,2^n]\nonumber
 \end{eqnarray}
\vspace{-2mm} and \vspace{-2mm}
\begin{equation} \label{H}
H=\delta_{2}[\overbrace{1,1,\cdots,1}^{2^{n-1}},\overbrace{2,2,\cdots,2}^{2^{n-1}}].\nonumber
\end{equation}
A straightforward computation shows that
$\mathcal{C}=\{\{i\}\}_{i=1}^{2^n}$, i.e. $|C|=1,\ \forall\ C\in
\mathcal{C}$. Thus the shift register is observable and
undecomposable with respect to outputs. To give an intuitive
explanation, we list the procedure for the case of
$n=3$ as follows . \vspace{-2mm}
$$
\begin{array}{rcl}
H&=&\delta_{2}[1,1,1,1,2,2,2,2],\nonumber\\
HL_i&=&\delta_{2}[1,1,2,2,1,1,2,2]\ \  (i=1,2),\nonumber \\
HL_iL_j&=&\delta_{2}[1,2,1,2,1,2,1,2]\ \ (i,j=1,2),\nonumber \\
HL_{i}L_{j}L_{k}&=&\!\!\!\left\{\begin{array}{l}
                                \!\!\!\delta_{2}[1,1,1,1,1,1,1,1],\ \ k=1,\\
                                \!\!\!\delta_{2}[2,2,2,2,2,2,2,2],\ \ k=2.\end{array}\right.
\end{array}
$$
It is easy to check that all the columns of $\mathcal{O}$ are distinct.
%
\par In order to display the procedure of finding a CC-PEVP,
we reconsider the BCN given in Example 1.\par {\bf Example 3.} For
BCN (\ref{eq25-00}), we have \vspace{-3mm} \vspace{0mm}
\begin{eqnarray*}
H &\!=& \!\delta_2[2,~1,~1,~1,~1,~1,~1,~2], \\
H\!L_1 &\!=& \!\delta_2[1,~2,~1,~2,~2,~1,~2,~1], \\
 H\!L_2^2\!=\!H\!L_2L_1\!=\!H\!L_1^2\!=\!H\!L_2 &\!=& \!\delta_2[1,~1,~1,~1,~1,~1,~1,~1], \\
H\!L_1L_2 &\!=& \!\delta_2[2,~2,~2,~2,~2,~2,~2,~2].
\end{eqnarray*}
From the observability matrix, we have
\begin{equation}
\label{eqC}
 \mathcal{C}=\{\{1,8\},\ \{2,4,5,7\},\ \{3,6\}\}.
\end{equation}
By $(\ref{eqC})$, we first let $S_1=\{1,8\}$, $S_2=\{3,6\}$. Since
$\mathcal{N}^2(S_1)=\{4,5\}$, we let $S_3=\{4,5\}$ and
$S_4=\{2,7\}$. It is easy to check that $\{S_l\}_{l=1}^4(|S_l|=2)$ is
the unique CC-PEVP. Thus, the system is decomposable with respect
to outputs. In fact, the partition obtained here is the same one as
shown in Fig.1.
\par
A contribution of this paper lies in that the regularity
assumption on the unobservable subspace proposed by Cheng, Li, \&
Qi (2010) is removed. To illustrate this point, we consider an example as
follows.   \par {\bf Example 3.} Consider BCN (\ref{eq4}) with
%
%
$L\!\!=\!\![L_{1}, L_{2}]\!\!\in\!\mathcal{L}_{8\times \!16}$,
$$L_{1}=\delta_{8}[6,8,1,8,7,8,6,8], \ L_{2} =\delta_{8}[6,  8,  7,  8,  1,  8,  6,  8]$$
and $H =\delta_{2}[1,  1,  2,  1,  2,  1,  1,  1].$ A
straightforward computation shows that
\begin{equation}
\mathcal{O}=\left[\begin{array}{cccccccc}
 1&1&2&1&2&1&1& 1\\
1&1&1&1&1&1&1& 1
\end{array}\right].
\end{equation}
Thus
$$
 \mathcal{C}=\{\{1, 2,  4, 6, 7, 8\},\ \{3, 5\}\}.
$$
We try to search a CC-PEVP $\mathcal{S}$ from $\mathcal{C}$. Let
$S_1=\{3,5\}$. It is easy to check that
$\mathcal{N}^1(S_1)=\mathcal{N}^2(C_2)=\{1,7\}.$ So we let
$S_2=\{1, 7\}$. Assume $S_3=\{2, 6\}$ and $S_4=\{4, 8\}$. Then it
is easily seen that
$$
\mathcal{N}^1(S_1)=\mathcal{N}^2(S_1)=\{1,7\}=S_2,
$$
$$
\mathcal{N}^1(S_2)=\mathcal{N}^2(S_2)=\{6\}\subset S_3,
$$
$$
\mathcal{N}^1(S_3)=\mathcal{N}^2(S_3)=\{8\}\subset S_4,
$$
$$
\mathcal{N}^1(S_4)=\mathcal{N}^2(S_4)=\{8\}\subset S_4.
$$
Furthermore, by the color of vertices, the vertices in $S_l$ have
the same output. Thus the system is decomposable with respect to outputs of
the maximum order 1. From $\{S_i\}_{i=1}^4$, we get
$Q=\delta_{4}[2,3,1,4,1,3,2,4]$. Let
\begin{equation*}
(I_{4}\otimes\textbf{1}_{2}^{\mathrm{T}})T=Q,
\end{equation*}
which implies $T=\delta_{8}[3,  5,  1,  7,  2,  6,  4,  8 ]$. Thus
the logical coordinate transformation $z=Tx$ realizes the maximum
decomposition with respect to outputs. It is worth
noting that the CC-PEVP is not unique. Another one is $\{\tilde
S_i\}_{i=1}^4$ with
$$
\tilde S_1=\{3,5\},\ \tilde S_2=\{1, 7\},\tilde S_3=\{2, 4\},\ \tilde S_4=\{6, 8\},
$$
which results a coordinate transformation matrix $\tilde T=\delta_{8}[3,  5,  1,  6,  2,  7,  4,  8 ].$ A straightforward computation shows that
$$R=\frac{1}{4}(\textbf{1}_{4}^{\mathrm{T}}\otimes
I_{2})T\tilde T^{\mathrm{T}}(\textbf{1}_{4}\otimes
I_{2})=\frac{1}{4}\left[\begin{array}{cc}
3 &   1\\
1 &   3\end{array}\right].$$ Since $R$ is not a logical
matrix, by Remark 1, the largest unobservable subspace is not regular.
\par

\section{Conclusions}

We have studied the decomposition with respect to outputs for BCNs, which is a generalization of the observability
decomposition of the traditional linear control theory. Our analysis relies on some equivalent algebraic and
graphical conditions for the decomposability with respect to outputs. It has been revealed that a BCN is decomposable with respect to outputs if and only if it has a CC-PEVP. With the observability matrix, an effective approach has been proposed for searching a CC-PEVP. In our future work, the Kalman decomposition without the regularity assumptions will be considered.

\end{document}